\documentclass[11pt]{article}
\usepackage[a4paper,margin=25mm,headheight=14pt]{geometry}
\usepackage[T1]{fontenc}
\usepackage{lmodern,amsmath,amssymb,amsthm,mathtools,microtype}
\usepackage{booktabs,array,longtable,needspace}
\usepackage{tikz}
\usepackage{subfig}
\usetikzlibrary{arrows.meta,positioning,calc}
\usepackage{fancyhdr}
\usepackage[hidelinks]{hyperref}
\usepackage[nameinlink,noabbrev]{cleveref}
\usepackage{comment}
\usepackage{changepage}

\newtheorem{theorem}{Theorem}[section]
\newtheorem{lemma}[theorem]{Lemma}
\newtheorem{proposition}[theorem]{Proposition}
\newtheorem{observation}[theorem]{Observation}
\newtheorem{corollary}[theorem]{Corollary}
\theoremstyle{definition}
\newtheorem{definition}[theorem]{Definition}
\theoremstyle{remark}

\crefname{lemma}{lemma}{lemmas}
\Crefname{lemma}{Lemma}{Lemmas}
\crefname{corollary}{corollary}{corollaries}
\Crefname{corollary}{Corollary}{Corollaries}
\crefname{proposition}{proposition}{propositions}
\Crefname{proposition}{Proposition}{Propositions}
\crefname{definition}{definition}{definitions}
\Crefname{definition}{Definition}{Definitions}
\crefname{remark}{remark}{remarks}
\Crefname{remark}{Remark}{Remarks}

\newtheoremstyle{boldcase}% name of the style to be used
{\topsep}% measure of space to leave above the theorem. E.g.: 3pt
{\topsep}% measure of space to leave below the theorem. E.g.: 3pt
{\normalfont}% name of font to use in the body of the theorem
{0pt}% measure of space to indent
{\bfseries}% name of head font
{.}% punctuation between head and body
{ }% space after theorem head; " " = normal interword space
{\thmname{#1}\thmnumber{ #2}:\thmnote{ #3}}% Manual spec for case layout to keep it bold

\theoremstyle{boldcase}
\newcounter{case}
\newcounter{subcase}[case]

\renewcommand{\thecase}{\arabic{case}}
\renewcommand{\thesubcase}{\thecase.\arabic{subcase}}

\newenvironment{case}{%
	\refstepcounter{case}%
	\par\medskip
	\noindent\textbf{Case \thecase.}\enspace
	\ignorespaces
}{%
	\par\medskip
}

\newenvironment{subcase}{%
	\refstepcounter{subcase}%
	\begin{adjustwidth}{1.5em}{0pt}%
		\noindent\textbf{Subcase \thesubcase.}\enspace
		\ignorespaces
	}{%
		\par
	\end{adjustwidth}%
}

\newcommand{\GD}{G^{\mathrm D}}

\title{Intersections of Odd Polygons Revisited}
\author{Eyal Ackerman\thanks{Department of Mathematics, Physics, and Computer Science, Oranim College, Tivon 36006, Israel. eyal\_a@oranim.ac.il}}

\begin{document}
\maketitle

\begin{abstract}
It is shown that the maximum number of intersection points between the boundaries of a simple $m$-gon and a simple $n$-gon is $mn-(m+n)+3$ for any odd integers $m,n \ge 3$.
The previous best upper bound was $mn-(m+n)+C$, for some huge constant $C$.
\end{abstract}

%\setcounter{tocdepth}{1}
%\tableofcontents
%\newpage

\section{Introduction}
\begin{figure}[t]
	\centering
	\includegraphics[scale=0.99]{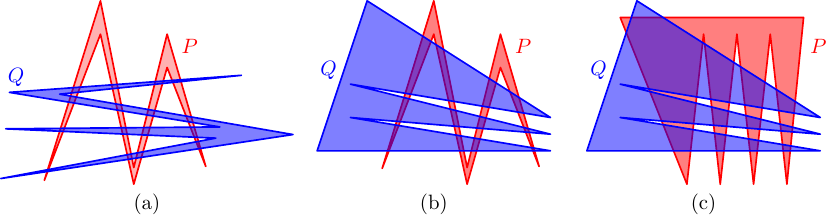}
	\caption{%
		(a)~Optimal construction for $m=n= 8$, with $8\cdot 8=64$ intersections.
		(b)~Optimal construction for $m= 8$, $n=7$, with $8\cdot 7-8 =48$ intersections.
		(c) Lower-bound construction for $m=9$, $n=7$, with
		$9\cdot 7 -(9+7) + 3 = 50$ intersections.
	}
	\label{fig:odd-odd-example}
\end{figure}

Let $P$ and $Q$ be two simple polygons with $m$ and $n$ sides, respectively.
We are interested in the very basic question of determining the maximum number of distinct intersection points between the boundaries of $P$ and $Q$.
We assume that no two sides of the polygons overlap, therefore, this number is finite and is upper-bounded by $mn$.
It is easy to see that this bound can be attained when $m$ and $n$ are both even, see Figure~\ref{fig:odd-odd-example}(a) for an example.
Observe also that if $n$ is odd, then every side of $P$ has at most $n-1$ distinct intersection points with sides of $Q$, therefore there are at most $mn-m$ intersections between sides of $P$ and $Q$.
It is not hard to construct an $m$-gon and an $n$-gon admitting that many intersections for every even $m$ and odd $n$, see Figure~\ref{fig:odd-odd-example}(b) for an example.
In contrast, determining the maximum number of intersections between the boundaries of two odd size polygons is much more challenging.

This problem was first studied by Dillencourt, Mount, and Saalfeld~\cite{DMS93}. 
They showed a construction with $mn-(m+n)+3$ intersections (see Figure~\ref{fig:odd-odd-example}(c) for an example) and claimed that this bound is tight, however, their proof contained an unrecoverable error.
A weaker upper bound of $mn-(m+\lceil n/6 \rceil)$ for $m \ge n$ was proved by \v{C}ern\'y, K\'ara, Kr\'al', Podbrdsk\'y, Sot\'akov\'a, and \v{S}\'amal~\cite{czechs}.
This was improved to $mn-(m+n)+C$, where $C$ is a (huge) constant, in joint work Keszegh and Rote~\cite{AKR22}.
Recently, the author revisited this work with ChatGPT, asking it whether that proof is sound and whether the constant $C$ can be decreased, hopefully, down to $3$.
The outcome of the first question was a small gap revealed in a certain case-analysis in the proof of the main theorem where one subcase was overlooked (see Appendix~\ref{app:gap}). 
Still, it seems that by a more careful and finer analysis this gap can be resolved.
As for the second question, using tools from~\cite{AKR22}, ChatGPT indeed managed to decrease the constant $C$: first, to $11$ while also providing a simpler proof than the one in~\cite{AKR22}. Then, with some extra effort and case-analysis, down to $3$, thus determining the exact maximum possible number of intersections between two odd size simple polygons.

\begin{theorem}\label{thm:main}
	For any odd integers $m,n \ge 3$, the maximum number of intersection points between the boundaries of a simple $m$-gon and a simple $n$-gon is $mn-(m+n)+3$.
\end{theorem}

While preparing this manuscript the author was contacted by Ju Yong Seo who published a different proof of Theorem~\ref{thm:main} obtained using a different LLM~\cite{Claude}.
In fact, the result in~\cite{Claude} is stronger, since it states that the so-called \emph{disjointness graph} of the two polygons has at most three connected components and this immediately implies Theorem~\ref{thm:main} (see Section~\ref{sec:tools}).
However, the proof in~\cite{Claude} is computer-assisted (in the pre-LLM era sense), as it boils down to verifying a few hundred cases using a SAT solver, whereas the proof presented here can be verified by hand.
Furthermore, we present a simple proof (relying on~\cite{AKR22}) of the slightly weaker bound $mn-(m+n)+11$.

\paragraph{Organization.}
We first establish terminology and some tools in Section~\ref{sec:tools}.
Using these tools we present a simple proof of the weaker bound in Section~\ref{sec:weaker}.
Section~\ref{sec:more-tools} introduces further definitions and tools that are required for proving the optimal bound which is then proved in Section~\ref{sec:optimal}.
The small gap in the proof in~\cite{AKR22} is described in Appendix~\ref{app:gap}.

\paragraph{AI Disclosure.}
ChatGPT-6 Astra was used for proving the new result and for preparing this manuscript.
The author verified the proofs, simplified some of them, edited the paper considerably to improve the presentation, and takes full responsibility for its content.

\section{Assumptions, terminology, and tools}
\label{sec:tools}

Throughout this paper $P$ denotes a simple $m$-gon and $Q$ denotes a simple $n$-gon, where $m$ and $n$ are odd integers.
We assume that each side of $P$ intersects each side of $Q$ in at most one point.
Furthermore, it follows from~\cite{AKR22} that we may assume, without loss of generality, that each intersection point of two sides lies in the interior of both sides (that is, it is a crossing point), no three vertices of $P$ and $Q$ combined are collinear, and that no two sides of $P$ and $Q$ combined are parallel.
Denote by $X$ the number of intersection points between a side of $P$ and a side of $Q$.

\paragraph{The disjointness graph $\GD$.}
As in \cite{AKR22,czechs}, we define the \emph{disjointness graph}
of $P$ and $Q$.
It is a bipartite graph, $\GD$, whose vertex set corresponds to the sides of $P$ and $Q$
and its edge set consists of pairs of disjoint sides, one from $P$ and the other from $Q$. See Figure~\ref{fig:H} for an example.
By abuse of notation we interchange sides of the polygons with their corresponding vertices in $\GD$.

\medskip
Let $c(\GD)$ be the number of connected components of $\GD$. % and let $CC(x)$ denote the connected component of $\GD$ that contains the polygon-side $x$.
Since every edge in $\GD$ corresponds to a non-intersecting pair of sides, it follows that $X = mn-|E(\GD)|$.
Note that $|V(\GD)|=m+n$, and therefore $|E(\GD)| \ge m+n-c(\GD)$.
Thus, we have:
\begin{equation}\label{eq:cc-bound}
X \le mn-(m+n)+c(\GD).
\end{equation}

The upper bound in~\cite{AKR22} follows by showing that $c(\GD) < 2^{2^{70}}$.
The current improvement of this bound is obtained by showing that many connected components imply many additional edges within connected components.
Clearly, once the number of edges in a graph exceeds the number of vertices minus one, any additional edge within a connected component implies a new cycle.
It is therefore convenient to consider the \emph{cycle rank} of $\GD$ and its connected components.
For a graph $G$ the \emph{cycle rank} of $G$ (also known as \emph{circuit rank} or \emph{cyclomatic number}) is the minimum number of edges one has to remove from $G$ to make it acyclic (that is, a forest).
We recall the following facts about the cycle rank of a graph, which we denote here by $\mu(G)$.

\begin{proposition}\label{prop:cycle-rank}
	Let $G$ be a graph, possibly with parallel edges.
	Then $\mu(G)=|E(G)|-|V(G)|+c(G)$, where $c(G)$ is the number of connected components of $G$.
	Furthermore, $\mu(G)$ is non-negative, does not decrease when adding vertices or edges, and remains the same when contracting a non-loop edge. If $G_1,G_2$ have disjoint edge sets, then $\mu(G_1 \cup G_2) \ge \mu(G_1) + \mu(G_2)$.
\end{proposition}

\begin{proof}
	If $c(G)=1$, that is $G$ is connected, then clearly $\mu(G)=|E(G)|-|V(G)|+c(G)$ since one has to remove $|E(G)|-(|V(G)|-1)$ edges to get a spanning tree with $|V(G)|-1$ edges.
	If $c(G)>1$, then by applying this argument for every connected component and summing up we get the desired equality.

	By definition $\mu(G)$ is non-negative and can only increase when adding vertices or edges.
	Contracting a non-loop edge $e$ decreases both the number of vertices and the number of edges by one while leaving the number of connected components unchanged, therefore $\mu(G/e) = \mu(G)$.
	It is also easy to see that $\mu(G_1 \cup G_2) \ge \mu(G_1) + \mu(G_2)$, since breaking every cycle in $G_1 \cup G_2$ requires in particular breaking every cycle in $G_1$ and every cycle in $G_2$.
\end{proof}

Since $\mu(\GD) = |E(\GD)|-|V(\GD)| + c(\GD)$, we have:
\begin{equation}\label{eq:cycle-rank-bound}
	X = mn-(m+n)+c(\GD)-\mu(\GD).
\end{equation}

Therefore, we can get a better upper bound on $X$ by showing that $\GD$ has few connected components or that its cycle rank is high.
More precisely, Theorem~\ref{thm:main} will follow once we show that $c(\GD)-\mu(\GD) \le 3$.
In fact, it is enough to show that $c(\GD)-\mu(\GD) \le 4$ by the following observation.

\begin{proposition}\label{prop:odd-deficit}
	$c(\GD)-\mu(\GD)$ is odd.
\end{proposition}

\begin{proof}
	The number of crossings between two Jordan curves is even, thus $X$ is even.
	Since $mn$ is odd and $m+n$ is even, it follows from~(\ref{eq:cycle-rank-bound}) that $c(\GD)-\mu(\GD)$ is odd.
\end{proof}

\begin{figure}
	\centering
	\includegraphics[width=7cm]{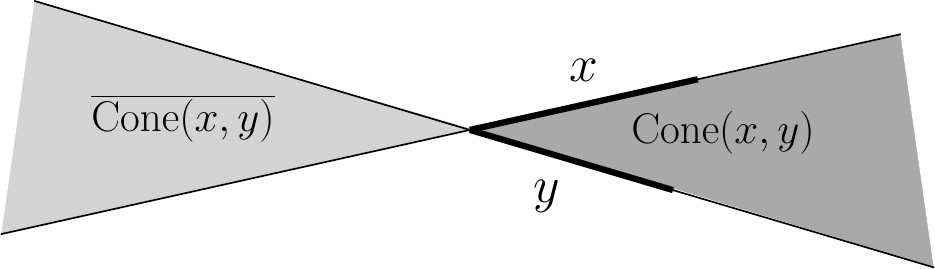}
	\caption{The cones defined by two consecutive polygon-sides $x$ and $y$.}
	\label{fig:cones}
\end{figure}

We proceed with some other useful notions and lemmas from~\cite{AKR22}.
Let $x$ and $y$ be two consecutive polygon-sides sharing a common endpoint $x\cap y$.
Denote by $\mbox{Cone}(x,y)$ the convex cone which is bounded by the two rays emanating from $x\cap y$ and containing $x$ and $y$, respectively, and let $\overline{\mbox{Cone}(x,y)}$ denote the opposite cone (see Figure~\ref{fig:cones}).

\begin{observation}\label{obs:out-of-cones}
	Let $x$ and $y$ be two consecutive polygon-sides and let $z$ be a polygon-side that crosses both $x$ and $y$.
	Then no endpoint of $z$ lies in $\mbox{Cone}(x,y) \cup \overline{\mbox{Cone}(x,y)}$.
\end{observation}

The following lemma gives rise to the main geometric structure used throughout the proof.

\begin{lemma}[{\cite[Lemma 5]{AKR22}}]\label{lem:gadget}
	Let $x$ and $x'$ be two consecutive sides of one polygon such that there is no side of the other polygon which is disjoint from both of them.
	Then there are two consecutive sides of the other polygon, $y$ and $y'$, such that $(x,y),(x',y') \in E(\GD)$ and $(x,y'),(x',y) \notin E(\GD)$.
	Furthermore, $x \cap x' \in \mbox{Cone}(y,y')$ or $y \cap y'\in \mbox{Cone}(x,x')$. 
\end{lemma}

We remark that the original statement of this lemma in~\cite{AKR22} assumes that $x$ and $x'$ belong to different connected components of $\GD$, however, its proof only uses the property stated here which is implied by the different components assumption.

\begin{figure}[t]
	\centering
	\includegraphics[width=16cm]{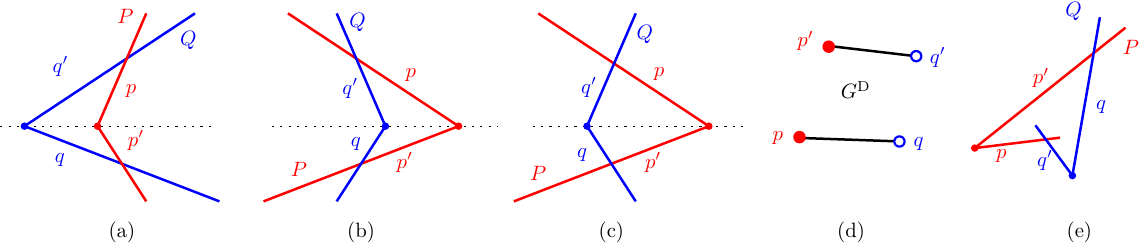}
	\caption{A gadget $(p,p',q,q')$.
		(a)~A pure $Q$-hooking gadget.
		(b)~A pure $P$-hooking gadget.
		(c)~A double gadget.
		(d)~The base edges of the gadget.
		(e)~Not a gadget since no hooking occurs.
	}
	\label{fig:hooked}
\end{figure}

\paragraph{Gadgets and hooking types.}
Let $q$ and $q'$ be two consecutive sides of $Q$ and let $p$ and $p'$ be two consecutive sides of $P$ such that $(q,p),(q',p') \in E(\GD)$, $(q,p'),(q',p) \notin E(\GD)$ and  $p \cap p' \in \mbox{Cone}(q,q')$ or $q \cap q'\in \mbox{Cone}(p,p')$ (see Figure~\ref{fig:hooked} for examples).
The two pairs of consecutive edges are called \emph{associated pairs} in~\cite{AKR22}.
Here we refer to them as a \emph{gadget} and say that $(p,q)$ and $(p',q')$ are its \emph{base} edges, or simply \emph{bases}.
As a convention, we denote a gadget by a $4$-tuple $(p,p',q,q')$ such that $p,p'$ are its $P$-pair, $q,q'$ are its $Q$-pair, and $(p,q),(p',q') \in E(\GD)$ are its base edges.
Note that by definition we get the same set of gadgets if we interchange $P$ and $Q$.

If the common endpoint of the $Q$-pair (resp., $P$-pair) of a gadget lies in the cone of the $P$-pair (resp., $Q$-pair), then we say that the  gadget is \emph{$P$-hooking} (resp., \emph{$Q$-hooking}).%\footnote{The notion of hooking is similar in~\cite{AKR22}, but not identical.} 
Note that a gadget may be both $P$-hooking and $Q$-hooking and that a pair of consecutive sides may participate in several gadgets (of similar or different hooking types).
If exactly one hooking occurs in a gadget, then we say that this gadget is \emph{pure} and otherwise if the two hookings occur the gadget is \emph{double}.
See Figure~\ref{fig:hooked} for illustrations.

\begin{observation}\label{obs:cones}
	Let $g=(p,p',q,q')$ be a gadget.\\
	If $g$ is pure $Q$-hooking, then $p \cap p' \in \mbox{Cone}(q,q')$ and $q \cap q' \in \overline{\mbox{Cone}(p,p')}$.\\
	If $g$ is pure $P$-hooking, then $q \cap q' \in \mbox{Cone}(p,p')$ and $p \cap p' \in \overline{\mbox{Cone}(q,q')}$.\\
	If $g$ is double hooking, then $q \cap q' \in \mbox{Cone}(p,p')$ and $p \cap p' \in \mbox{Cone}(q,q')$.
\end{observation}

\begin{corollary}\label{cor:no-isolated}
	$\GD$ contains no isolated vertices.
\end{corollary}

\begin{proof}
	Suppose for contradiction that there is a polygon-side $x$ which is an isolated vertex in $\GD$. Let $x'$ be a side of the same polygon that is consecutive to $x$.
	Then it follows from Lemma~\ref{lem:gadget} that there is a side $y$ of the other polygon such that $(x,y) \in E(\GD)$, a contradiction.
\end{proof}

The following lemma is crucial for the proof. 

\begin{lemma}[{\cite[Lemma 7]{AKR22}}]\label{lem:eight}
	Suppose that two gadgets use eight distinct sides and induce a subgraph of $\GD$ which is a matching consisting of their four base edges. 
	Then both gadgets are pure and of different hooking types.
\end{lemma}

\paragraph{Remark.} We note that the original statement of this lemma refers to the case where the four base edges belong to different connected components of $\GD$, however, the proof actually uses the weaker assumption stated here.

\section{A weaker upper bound $X \le  mn-(m+n)+11$}
\label{sec:weaker}

Recall that it is enough to show that $c(\GD)-\mu(\GD) \le 11$.
In particular we may assume that $c(\GD) \ge 2$.

\begin{figure}[t]
	\centering
	\includegraphics[scale=1]{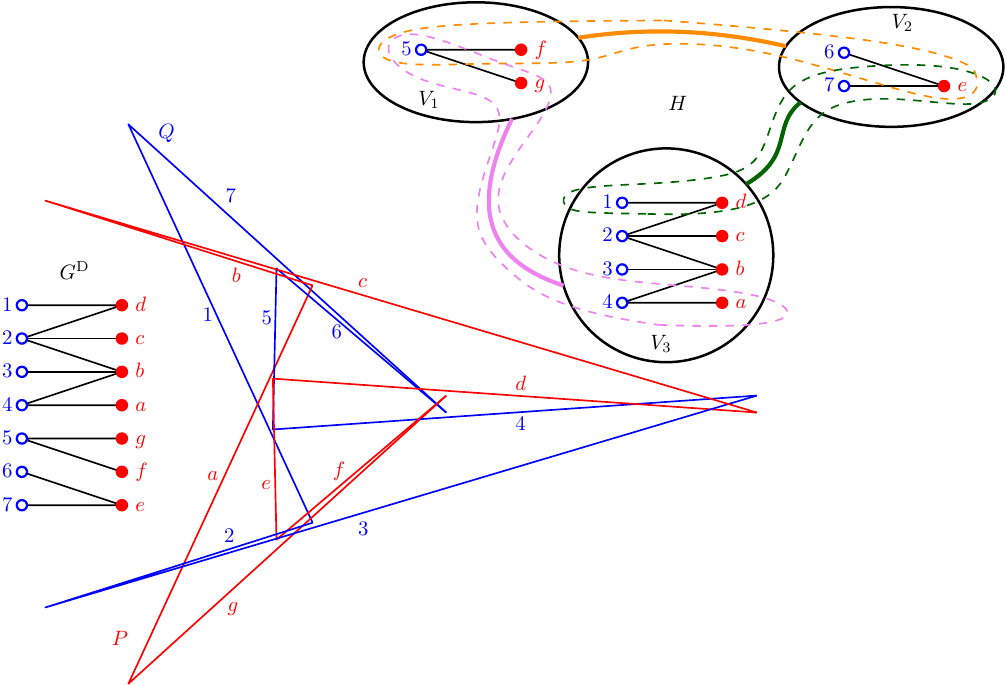}
	\caption{Two heptagons $P$ and $Q$, their disjointness graph $\GD$, and their transitions graph $H$. The edges of $H$ are associated with three gadgets: the $P$-pure gadget $(a,g,4,5)$, the $Q$-pure gadget $(d,e,1,7)$ and the double gadget $(e,f,6,5)$. 
	%$H_P=(\{V_1,V_2,V_3\},\{V_1V_3,V_1V_2\})$ and $H_Q=(\{V_1,V_2,V_3\},\{V_2V_3,V_1V_2\})$.
	}
	\label{fig:H}
\end{figure}

\paragraph{The transitions graph $H$.}
Define a graph $H$ as follows:
the node\footnote{We will refer to the vertices of $H$ as `nodes' in order to distinguish them from the vertices of $\GD$.} set of $H$ consists of the connected components of $\GD$ and its edge set consists of pairs of distinct connected components $U$ and $V$ such that there are two consecutive $Q$-sides $q$ and $q'$ such that $q\in U$ and $q' \in V$.
By Lemma~\ref{lem:gadget} there is at least one gadget whose $Q$-pair is $q,q'$. 
We say that any such gadget and $UV$ are \emph{associated} and note that 
the gadget has one of its base edges in $U$ and its other base edge in $V$.
See Figure~\ref{fig:H} for an example.

\paragraph{Remarks.}
It easily follows from its definition that $H$ is connected. 
Observe that a gadget is associated with at most one edge of $H$, but several gadgets may be associated with the same $H$-edge.
The fact that an $H$-edge corresponds to a gadget and the symmetric roles that $P$ and $Q$ play with respect to gadgets imply that defining the transitions graph according to $P$ would result in the same graph $H$.

%Let $H_Q$ (resp., $H_P$) be the subgraph of edges that are associated with a $Q$-hooking (resp., $P$-hooking) gadget.
%Note that $H = H_Q \cup H_P$ and that an edge may belong to both subgraphs.
%Clearly, every gadget is associated with exactly one edge.

%We say that two edges in a graph are \emph{near} if they share at least one endpoint or there is an edge connecting an endpoint of one edge to an endpoint of the other edge.

%\begin{observation}\label{obs:degree-2}
%	Let $UV_1$ and $UV_2$ be distinct edges of $H$ and let $g_1=(p,p_1,q_1,q'_1)$ and $g_2=(p,p_2,q_2,q'_2)$ be gadgets associated with them, respectively, such that $(p,q_1)$ and $(p,q_2)$ are their bases in $U$. Then $p_1 \ne p_2$. 
%\end{observation}

%\begin{proof}
%	The claim holds since $p_1\in V_1$ and $p_2\in V_2$.
%\end{proof}

%Since every polygon-side has two neighboring sides that we conclude the following.

%\begin{corollary}\label{cor:degree-2}
%	In the subgraph of $\GD$ that consists of the base edges of gadgets associated with edges of $H$ the degree of every vertex is at most two.
%\end{corollary}

%\begin{proof}
%	Suppose without loss of generality that a $P$-side $p$ in a connected component $U$ belongs to three distinct bases $(p,q_1)$, $(p,q_2)$ and $(p,q_3)$. 
%	Let $(p,p_i,q_i,q'_i)$ be the corresponding gadget and let $UV_i$ be the associated $H$-edge, for $i=1,2,3$.
%	Then $p_i$ is a consecutive side of $p$ and belongs to $V_i$, for $i=1,2,3$. 
%	But $p$ has only two consecutive sides.
%\end{proof}

\begin{lemma}\label{lem:same-type-not-disjoint}
	There are no two node-disjoint edges of $H$ associated with two gadgets that are both $Q$-hooking or both $P$-hooking.
\end{lemma}

\begin{proof}
	Suppose that $H$ has two disjoint edges associated with same-type gadgets $g_1$ and $g_2$.
	Then the four bases of $g_1$ and $g_2$ lie in four distinct connected components, thus, their eight vertices induce a matching in $\GD$. But then by Lemma~\ref{lem:eight} the two gadgets must have different types.
\end{proof}

\begin{lemma}\label{lem:near-bases}
	If two gadgets that join a component $U$ to two distinct other components are both $Q$-hooking or both $P$-hooking, then their bases in $U$ are connected by an edge of $\GD$ or share a vertex.
\end{lemma}

\begin{proof}
	Suppose that there are two gadgets sharing the same type $g_1=(p_1,p'_1,q_1,q'_1)$ and $g_2=(p_2,p'_2,q_2,q'_2)$ and three distinct connected components $U,L_1,L_2$ such that $(p_i,q_i)$ lies in $U$ and $(p'_i,q'_i)$ lies in $L_i$, for $i=1,2$
	(see Figure~\ref{fig:forced-disjoint} for an illustration).
	It is enough to consider the case that $p_1 \ne p_2$ and $q_1 \ne q_2$, for otherwise the claim trivially holds.
	Consider the three subsets $\{p_1,p_2,q_1,q_2\} \subseteq U$, $\{p'_1,q'_1\} \subseteq L_1$ and $\{p'_2,q'_2\} \subseteq L_2$, and observe that there cannot be edges in $\GD$ between vertices from different subsets, since they belong to different connected components.
	If $(p_1,q_2),(p_2,q_1) \notin E(\GD)$, then the subgraph induced by those eight vertices is a matching and by Lemma~\ref{lem:eight} the two gadgets cannot share a hooking type, a contradiction.
\begin{figure}[t]
	\centering
	\includegraphics[width=7cm]{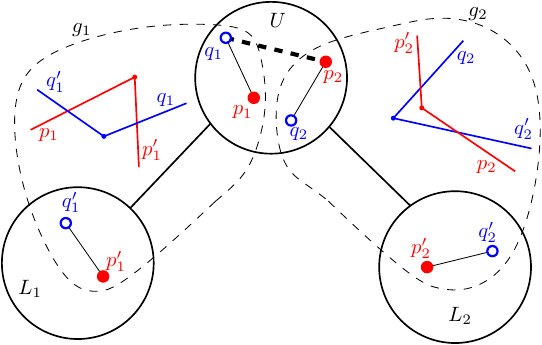}
	\caption{An illustration for the proof of Lemma~\ref{lem:near-bases}: the gadget $g_i$ consists of the $Q$-pair $q_i$,$q'_i$ and the $P$-pair $p_i,p'_i$, its base edges $(p_i,q_i)$ and $(p'_i,q'_i)$ are in $U$ and $L_i$, respectively. It follows that $(p_1,q_1)$ and $(p_2,q_2)$ share a vertex or are connected by an edge of $\GD$.}
	\label{fig:forced-disjoint}
\end{figure}
\end{proof}

The following is an important consequence of Lemma~\ref{lem:near-bases}.

\begin{corollary}\label{cor:matching}
	Let $U$ be a connected component containing $k\ge 1$ pairwise vertex-disjoint bases belonging to gadgets sharing the same hooking type where each of these gadgets is associated with an $H$-edge connecting $U$ to a distinct connected component.
	Then the vertices of the $k$ bases induce a subgraph of $U$ with $2k$ vertices and at least $k+\binom{k}{2}$ edges.
	In particular, $\mu(U) \ge k+\binom{k}{2}-(2k-1) = \binom{k-1}{2}$.
\end{corollary}

In order to use Corollary~\ref{cor:matching} we first consider a spanning tree $T$ of $H$.
Clearly, $T$ consists of exactly $c(\GD)-1$ edges.
For every edge $UV$ of $T$, we select one gadget associated with $UV$, and color $UV$ with color $Q$ if this gadget is $Q$-hooking, and with color $P$ otherwise (if it is $P$-pure).
Let $T_Q$ and $T_P$ be the subgraphs consisting of the edges colored $Q$ and $P$, respectively.
Recall that by Lemma~\ref{lem:same-type-not-disjoint} any two edges of $H$ associated with same-type gadgets share a common node of $H$.
It follows that apart from isolated nodes, each of $T_Q$ and $T_P$ is a star.\footnote{Had we stayed with $H$ instead of switching to $T$, we would have also had to consider the case of a triangle, which is the other possible graph without node-disjoint edges. We will handle this case in the proof of the optimal bound in Section~\ref{sec:optimal}.}
Let $U_P$ and $U_Q$ be the centers of the stars $T_P$ and $T_Q$, respectively
(if one of these stars is empty, then we choose an arbitrary component as its center).
Note that it may happen that $U_P=U_Q$.

\paragraph{Obtaining disjoint base edges in $U_P$ and $U_Q$.}
For each edge $U_QV$ of $T_Q$, we mark its base edge in $U_Q$ that belongs to the gadget selected for $U_QV$.
Then we do the same for edges of $T_P$ and their corresponding base edges in $U_P$.
If a base edge belongs to two selected gadgets, then we mark it twice.
Considering the subgraph of $\GD$ that consists of these $c(\GD)-1$ marked edges (possibly with parallel edges), note that its maximum degree is at most two.
Indeed, suppose without loss of generality that a $P$-side $p$ in $U_P$ belongs to three marked bases $(p,q_1)$, $(p,q_2)$ and $(p,q_3)$. 
Let $(p,p_i,q_i,q'_i)$ be the corresponding gadget and let $UV_i$ be the associated $H$-edge, for $i=1,2,3$. 
Note that these are three distinct $H$-edges, since we selected one gadget per tree edge.
Then $p_i$ is a consecutive side of $p$ and belongs to $V_i$, for $i=1,2,3$. 
But $p$ has only two consecutive sides.

It follows that the subgraph of marked bases consists of paths and even cycles (possibly $2$-cycles).
Therefore, its edge set can be partitioned into two matchings by picking every other edge in a path or a cycle.
It follows that there is a matching $M \subseteq E(\GD)$ such that each of its (marked base) edges lies in $U_Q$ or $U_P$ and
\begin{equation}\label{eq:M}
	|M| \ge \frac{c(\GD)-1}{2}.
\end{equation}

\paragraph{Finishing the proof.}
Let $M_Q \subseteq M$ (resp., $M_P$) denote the base edges in $M$ that belong to selected gadgets that are associated to edges of $T_Q$ (resp., $T_P$).
Set $|M_Q|=k$ and $|M_P|=t$.
Considering the subgraph of $U_Q$ that is induced by the endpoints of the edges in $M_Q$, it follows from Corollary~\ref{cor:matching} that its cycle rank is at least $\binom{k-1}{2} \ge 2k-5$.\footnote{This inequality holds for every integer $k$. If $k=0$, then we cannot use Corollary~\ref{cor:matching}, but clearly in this case the cycle rank is zero which is greater than $2k-5$.}
Similarly, the cycle rank of the subgraph of $U_P$ induced by the endpoints of the $t$ base edges in $M_P$ is at least $2t-5$.
Since these two subgraphs are edge-disjoint and $k+t=|M|$, it follows from Proposition~\ref{prop:cycle-rank} and~(\ref{eq:M}) that $\mu(\GD) \ge 2k-5+2t-5 \ge c(\GD)-11$ and by (\ref{eq:cycle-rank-bound}) we get that $X \le mn-(m+n)+11$ as claimed.

\section{More tools for proving the optimal bound}
\label{sec:more-tools}

Before proving the optimal upper bound on $X$ in Section~\ref{sec:optimal}, we collect further useful claims.

\begin{proposition}\label{prop:3-arrangement}
	Let $\ell_1, \ell_2, \ell_3$ be three lines in general position and let $\cal A$ denote their arrangement.
	Let $\Delta$ be the single (triangular) bounded face of $\cal A$ and let $v_1v_2$ be an edge of $\Delta$.
	For $i=1,2$, let $s_i$ be a straight-line segment with one endpoint in the unbounded face whose intersection with $\Delta$ is $v_i$ and its other endpoint lies in the unbounded face whose intersection with $\Delta$ is its edge opposite to $v_i$.\footnote{Here we actually consider intersections of the closures of the faces.}
	If each of $s_1$ and $s_2$ crosses $v_1v_2$, then $s_1$ and $s_2$ cross each other.
\end{proposition}

\begin{proof}
	By applying an appropriate affine transformation we may assume without loss of generality that $\ell_1$ is the $x$-axis, $\ell_2$ is the $y$-axis, $\ell_3$ is the line $x+y=1$, $v_1=(0,0)$ and $v_2=(1,0)$.
	We denote the unbounded faces by $F_1,\ldots,F_6$ in clockwise manner, such that $F_1$ shares the edge $v_1v_2$ with $\Delta$.
	We may assume that the endpoints of $s_1$ lie in $F_2$ and $F_5$ and the endpoints of $s_2$ lie in $F_3$ and $F_6$.
	See Figure~\ref{fig:3-arrangement}.
	\begin{figure}
		\centering
		\includegraphics[width=8cm]{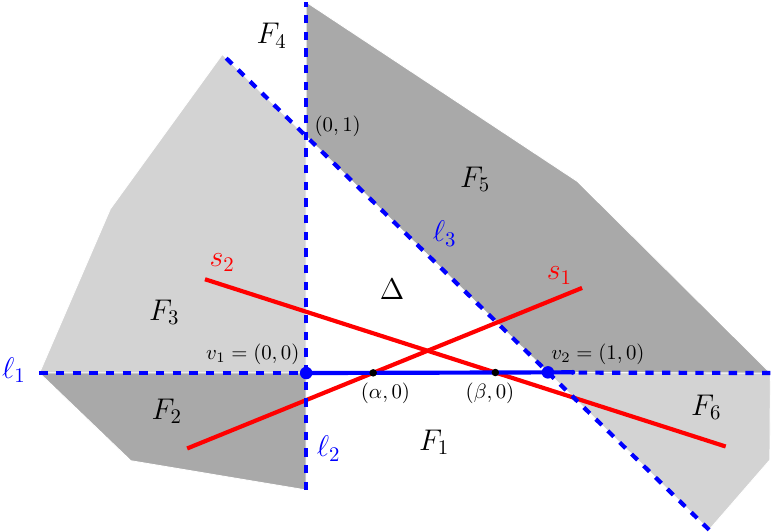}
		\caption{An illustration for the proof of Proposition~\ref{prop:3-arrangement}: If $s_1$ crosses $v_1v_2$ and its endpoints are in $F_2$ and $F_5$ and $s_2$ crosses $v_1v_2$ and its endpoints are in $F_3$ and $F_6$, then $s_1$ and $s_2$ cross.}
		\label{fig:3-arrangement}
	\end{figure}
	Following $s_1$ from its endpoint in $F_2$ it enters $\Delta$ at a crossing point with $v_1v_2$, denote it by $(\alpha,0)$, $0<\alpha<1$.
	Similarly, following $s_2$ from its endpoint in $F_6$ it enters $\Delta$ at a crossing point with $v_1v_2$, denote it by $(\beta,0)$, $0<\beta<1$.
	If $\alpha=\beta$, then $s_1$ and $s_2$ cross at $(\alpha,0)$.
	If $\beta < \alpha$, then $s_1$ and $s_2$ must cross inside $F_1$, since their crossing points with its boundary alternate.
	Otherwise, if $\beta > \alpha$, then $s_1$ and $s_2$ must cross inside $\Delta$ since their crossing points with its boundary alternate.
\end{proof}

\begin{definition}[near edges]
	Let $G$ be a graph and let $e_1$ and $e_2$ be two distinct edges in $G$.
	We say that $e_1$ and $e_2$ are \emph{near}, if they share a vertex or there is an edge between an endpoint of $e_1$ and an endpoint of $e_2$.
\end{definition}

Thus, Lemma~\ref{lem:near-bases} shows that there are near bases provided certain assumptions hold.
The following will be useful later on.

\begin{proposition}\label{prop:tree-packing}
	In a tree $T$:\\
	(a)~For every non-empty set of pairwise near edges in $T$, there is an edge $(a,b)$ such that every edge in the set has $a$ or $b$ as an endpoint. In particular $T$ has no three pairwise near edges that form a matching; and\\
	(b)~if $e_1$ and $e_2$ are two non-near edges, then a non-empty set of edges that are near each of $e_1$ and $e_2$ forms a star.
\end{proposition}

\begin{proof}
	Let $S$ be a set of pairwise near edges in $T$.
	If every pair of edges shares a vertex, then $S$ is a star (it cannot be a triangle), and we can pick any edge in $S$ as $(a,b)$.
	Otherwise, let $e_1$ and $e_2$ be vertex-disjoint edges and let $(a,b)$ be the edge connecting them.
	Suppose that there is an edge $e_3 \in S$ with no endpoint in $\{a,b\}$.
	Considering the connected components of $T \setminus \{(a,b)\}$, $e_3$ belongs to one of them which does not contain at least one of $e_1$ and $e_2$, since they belong to different components.
	Assume without loss of generality that $e_3$ and $e_1$ belong to different components. Then a path connecting them in $T$ must go through $(a,b)$ but only one of them is incident to this edge. 
	Therefore $e_1$ and $e_3$ are not near.
		
	As for the second part of the claim, consider the unique path in $T$ connecting $e_1$ and $e_2$ and let $v$ be a vertex on this path which is not an endpoint of any of them (there must be such a vertex since $e_1$ and $e_2$ are not near).
	We claim that every edge which is near both $e_1$ and $e_2$ must have $v$ as an endpoint.
	Indeed, suppose that there is an edge $e$ which is near both $e_1$ and $e_2$ and $v$ is not one of its endpoints. Considering the connected components of $T \setminus \{v\}$, $e$ belongs to one of them which does not contain at least one of $e_1$ and $e_2$, since they belong to different components.
	Assume without loss of generality that $e$ and $e_1$ belong to different components. Then a path connecting them in $T$ must go through $v$ and therefore they are not near edges.
\end{proof}

By Lemma~\ref{lem:near-bases} two side-disjoint gadgets spanning three connected components must be of different pure types or force an edge between their base edges in their common connected component.
The following claim complements this by showing that if the gadgets are both pure, then they do not share sides.

\begin{proposition}\label{prop:sharing}
	Let $g_1=(p_1,p'_1,q_1,q'_1)$ and $g_2=(p_2,p'_2,q_2,q'_2)$ be two
	gadgets and let $U,L_1,L_2$ be three distinct connected components of $G^D$ such that $(p_i,q_i)$ lies in $U$ and $(p'_i,q'_i)$ lies in $L_i$, for $i=1,2$. Then:\\
	if $p_1=p_2$, then both gadgets are $Q$-hooking and at least one is also $P$-hooking;\\
	if $q_1=q_2$, then both gadgets are $P$-hooking and at least one is also $Q$-hooking; and\\
	if $p_1=p_2$ and $q_1=q_2$, then both gadgets are double.\\
	Consequently, if $g_1$ and $g_2$ are pure, then their bases in $U$ are vertex-disjoint.
\end{proposition}

\begin{proof}
	Suppose that $p_1=p_2=p$.
	Apply an affine transformation of the plane such that the triangle whose vertices are $p \cap p'_1$, $p \cap q'_1$ and $p'_1 \cap q_1$ is mapped to the triangle whose vertices are $(0,0)$, $(1,0)$ and $(0,1)$, respectively.
	Note that if a point $x$ lies in a convex cone $y$, then the image of $y$ is a convex cone and the image of $x$ lies in it.
	Therefore, intersections and hooking-types remain the same under this transformation.
	
	We first claim that no gadget is $P$-pure.
	Suppose without loss of generality that $g_1$ is $P$-pure.
	Since $g_1$ is $P$-hooking, $q_1 \cap q'_1 \in \mbox{Cone}(p,p'_1)$.
	Since $g_1$ is not $Q$-hooking $q_1 \cap q'_1$ must lie inside the triangle whose corners are $p \cap p'_1$, $p \cap q'_1$ and $p'_1 \cap q_1$, see Figure~\ref{fig:no-P-pure}.
	However, then $p \cap p'_2 \in \mbox{Cone}(q_1,q'_1)$ which implies that $p'_2$ cannot cross both $q_1$ and $q'_1$ by Observation~\ref{obs:out-of-cones}.
	Thus, $p'_2 \in U$ or $p'_2 \in L_1$, a contradiction.
	\begin{figure}[t]
		\centering
		\includegraphics[width=5cm]{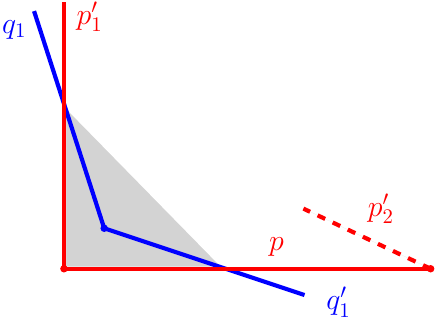}
		\caption{An illustration for the proof of Proposition~\ref{prop:sharing}: If $p_1=p_2=p$ and $g_1$ is $P$-pure, then $p'_2$ cannot cross both $q_1$ and $q'_1$.}
		\label{fig:no-P-pure}
	\end{figure}
	
	Therefore, both gadgets are $Q$-hooking and it remains to show that at least one of them is also $P$-hooking.
	Suppose that both of them are $Q$-pure.
	We proceed by considering whether $p'_2$ lies above or below the $x$-axis and whether it has a positive or negative slope.
	Consider the arrangement determined by the three lines through $p$, $p'_1$ and $p'_2
	$:
	it contains one bounded triangular face, denote it by $\Delta$, and six unbounded faces. 
	Denote these unbounded faces by $F_1,\ldots,F_6$ according to their clockwise order and such that $F_1$ has $p$ as one of its bounding edges.
	
	\begin{itemize}
		\item Suppose that $p'_2$ lies below the $x$-axis and has a negative slope, see Figure~\ref{fig:no-2-Q-pure4}.
		\begin{figure}[t]
			\centering
			\subfloat[If $p'_2$ has a negative slope, then $q_1 \cap q'_1 \in \mbox{Cone}(p,p'_2)$ and therefore $q'_1$ cannot cross both $p$ and $p'_2$.]{\includegraphics[width= 7cm]{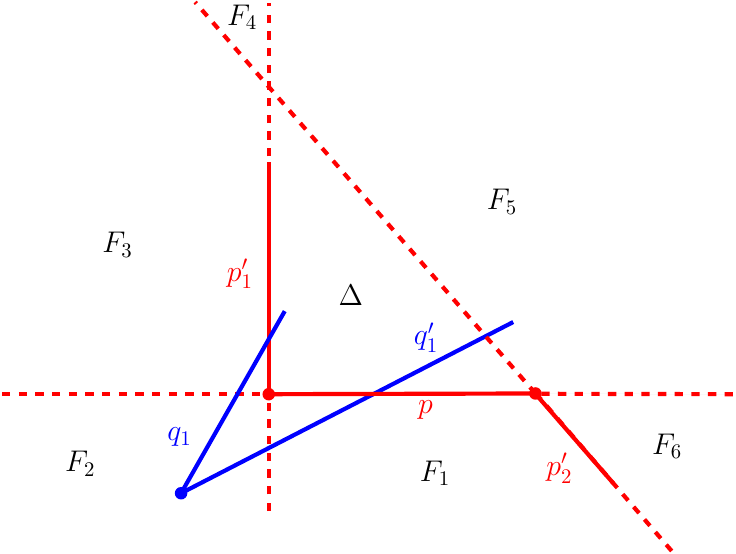}\label{fig:no-2-Q-pure4}}
			\hspace{5mm}
			\subfloat[If $p'_2$ has a positive slope, then $q_2 \cap q'_2 \in \mbox{Cone}(p,p'_1)$ and therefore $q'_2$ cannot cross both $p$ and $p'_1$.]{\includegraphics[width= 7cm]{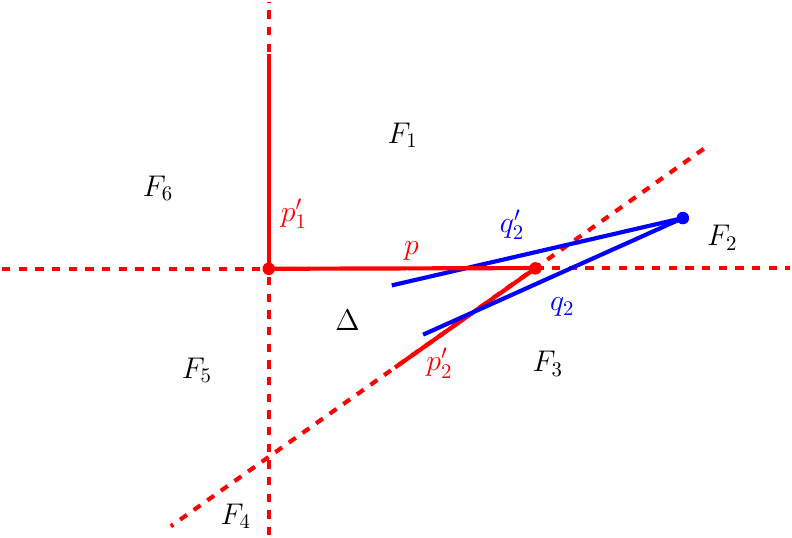}\label{fig:no-2-Q-pure3}}
			\caption{Illustrations for the proof of Proposition~\ref{prop:sharing}: $p'_2$ lies below the $x$-axis.}
			\label{fig:no-2-Q-pure-below}
		\end{figure}
		Then $\overline{\mbox{Cone}(p,p'_1)} = F_2$ and $q_1 \cap q'_1$ must lie in this face by Observation~\ref{obs:cones} since $g_1$ is $Q$-pure.
		However, then $q_1 \cap q'_1 \in \mbox{Cone}(p,p'_2)$ which implies by Observation~\ref{obs:out-of-cones} that $q'_1$ cannot cross both $p$ and $p'_2$, as it should.
		
		\item Suppose that $p'_2$ lies below the $x$-axis and has a positive slope, see Figure~\ref{fig:no-2-Q-pure3}.
		Then $\overline{\mbox{Cone}(p,p'_2)} = F_2$ and $q_2 \cap q'_2$ must lie in this face by Observation~\ref{obs:cones} since $g_2$ is $Q$-pure.
		However, then $q_2 \cap q'_2 \in \mbox{Cone}(p,p'_1)$ which implies that $q'_2$ cannot cross both $p$ and $p'_1$, as it should.
		
		\item Suppose that $p'_2$ lies above the $x$-axis and has a positive slope, see Figure~\ref{fig:no-2-Q-pure1}.
		\begin{figure}[t]
			\centering
			\subfloat[$p'_2$ has a positive slope.]{\includegraphics[width= 7cm]{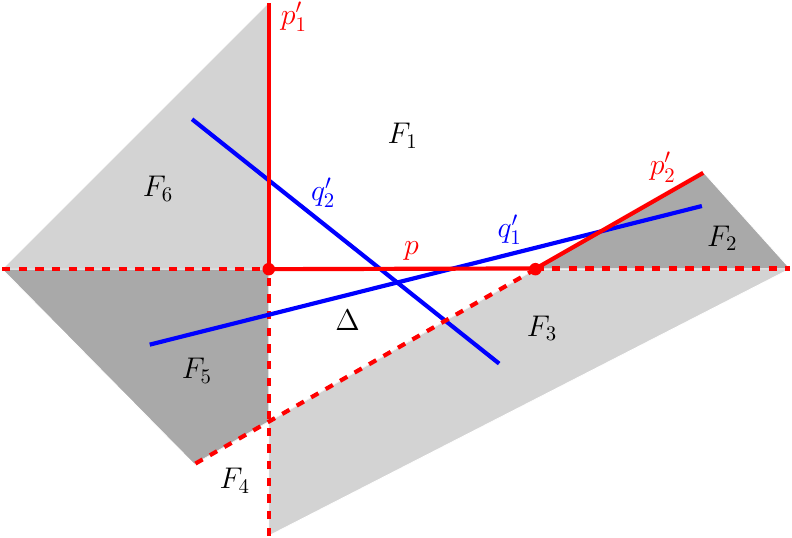}\label{fig:no-2-Q-pure1}}
			\hspace{5mm}
			\subfloat[$p'_2$ has a negative slope.]{\includegraphics[width= 6cm]{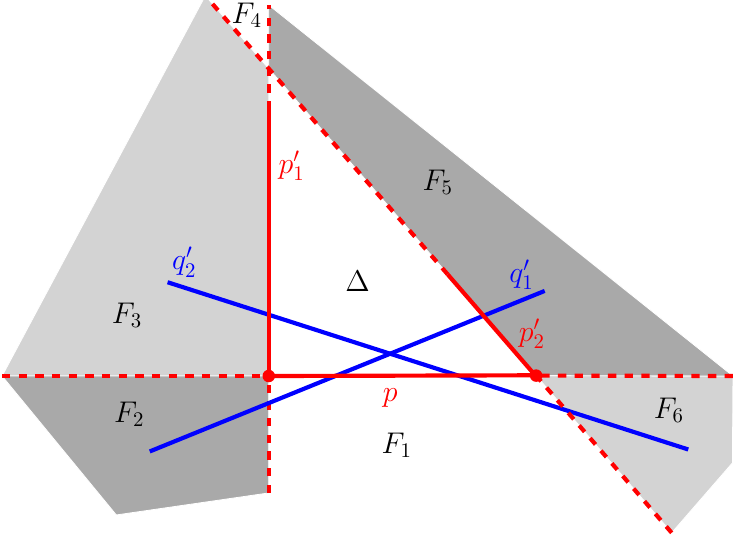}\label{fig:no-2-Q-pure2}}
			\caption{Illustrations for the proof of Proposition~\ref{prop:sharing}: $p'_2$ lies above the $x$-axis.}
			\label{fig:no-2-Q-pure-above}
		\end{figure}
		Since $g_1$ is $Q$-pure, $q_1 \cap q'_1 \in \overline{\mbox{Cone}(p,p'_1)}=F_4 \cup F_5$ by Observation~\ref{obs:cones}. Because $q'_1$ crosses $p$ and $p'_2$, it follows that $q'_1$ has one endpoint in $F_5$ and its other endpoint in $F_2$.
		Similarly, since $g_2$ is $Q$-pure, $q_2 \cap q'_2 \in \overline{\mbox{Cone}(p,p'_2)}=F_4 \cup F_3$. Because $q'_2$ crosses $p$ and $p'_1$, it follows that $q'_2$ has one endpoint in $F_3$	and its other endpoint in $F_6$.
		However, then $q'_1$ and $q'_2$ cross by Proposition~\ref{prop:3-arrangement}, which is impossible.
		
		\item Suppose that $p'_2$ lies above the $x$-axis and has a negative slope, see Figure~\ref{fig:no-2-Q-pure2}.
		Since $g_1$ is $Q$-pure, $q_1 \cap q'_1 \in \overline{\mbox{Cone}(p,p'_1)}=F_2$ by Observation~\ref{obs:cones}. Because $q'_1$ crosses $p$ and $p'_2$, it follows that $q'_1$ has its other endpoint in $F_5$.
		Similarly, since $g_2$ is $Q$-pure, $q_2 \cap q'_2 \in \overline{\mbox{Cone}(p,p'_2)}=F_6$. Because $q'_2$ crosses $p$ and $p'_1$, it follows that $q'_2$ has its other endpoint in $F_3$.
		However, then $q'_1$ and $q'_2$ cross by Proposition~\ref{prop:3-arrangement}, which is impossible.
	\end{itemize}
	
	This concludes the case that $p_1=p_2$.
	The case that $q_1=q_2$ is symmetric, and the case that $p_1=p_2$ and $q_1=q_2$ follows from the two previous cases.
	It also follows that if $g_1$ and $g_2$ are pure, then $p_1 \ne p_2$ and $q_1 \ne q_2$.
\end{proof}

Next, we establish further properties of the transition graph $H$.

\begin{proposition}\label{prop:cut}
For every partition of the nodes of $H$ into two disjoint subsets $A$ and $B$ such that $|A| \ge 2$ and $|B| \ge 1$ there are two gadgets crossing the cut with vertex-disjoint bases in $A$.
\end{proposition}

\begin{proof}
Observe first that there cannot be just one $Q$-side in a component from $A$ that is part of every gadget crossing the cut.
Indeed, suppose that there is such a side $q$ and let $V_1 \in A$ be the connected component of $\GD$ that contains $q$.
By deleting $q$ from $Q$, we get a path.
Therefore, if we construct $H$ using this path instead of $Q$, then we would still get a connected graph, apart from possibly $V_1$ which might be an isolated node.
Since $|A| \ge 2$, it follows that there are a side in $A$ and a side in $B$ that are consecutive along the path and therefore by Lemma~\ref{lem:gadget} belong to a gadget crossing the cut.  

Similarly, there cannot be just one $P$-side in a component from $A$ that is part of every gadget crossing the cut, since constructing the transitions graph with respect to $P$ results in the same graph.
Considering the $H$-edges crossing the cut and base edges in components from $A$ that belong to corresponding gadgets, we have a bipartite graph such that each of its classes (sides of $Q$ and $P$) contains at least two vertices and every vertex has degree at least one.
Such a bipartite graph has a matching of size two, that is, two vertex-disjoint bases in $A$.
\end{proof}

\begin{definition}
Let $U$ be a node of $H$ and let $S$ be a set of gadgets associated with $H$-edges $E(S) \subseteq E(H)$ adjacent to $U$.
We say that $S$ is \emph{admissible} (with respect to $U$) if each edge $e \in E(S)$ is associated with at most two gadgets in $S$ and if $e$ is associated with exactly two gadgets in $S$, then their bases in $U$ are vertex-disjoint.
\end{definition}

\begin{proposition}\label{prop:2-capacities}
Let $U$ be a node in $H$ and let $S$ be an admissible set of gadgets with respect to $U$.
Then $|S|\le 2\mu(U)+7$.
Furthermore, if all gadgets in $S$ are pure of the same type, then $|S| \le \mu(U)+3$.
\end{proposition}

\begin{proof}
	%We use induction on $\mu(U)$, assuming that $U$ satisfies all the properties established so far, regardless of its relation to intersecting polygons.
	Suppose first that $\mu(U)=0$, that is, $U$ is a tree.
	Let $S_Q$ denote the $Q$-hooking gadgets in $S$ and let $B_Q$ denote their bases in $U$, where a base may appear with multiplicity in $B_Q$ if it belongs to several gadgets in $S_Q$. 
	Let $d_Q$ denote the maximum vertex degree in the subgraph of $U$ consisting of the base edges in $B_Q$.
	Then $d_Q \le 2$: 
	Firstly, by admissibility a vertex cannot belong to two bases of two gadgets on the same $H$-edge.
	Secondly, if a vertex (polygon-side) belongs to three bases of gadgets on different $H$-edges, then a polygon-side incident to it must belong to the other nodes (components) of these edges. But that is impossible since a polygon-side has just two incident sides. 
	Now, if all the gadgets in $S_Q$ are $Q$-pure, then it follows from admissibility and Proposition~\ref{prop:sharing} that $B_Q$ is a matching, that is, $d_Q\le 1$. If the bases in $B_Q$ are pairwise near, then by Proposition~\ref{prop:tree-packing} there is an edge in $U$ that is incident to all of them. This implies that $|S_Q|=|B_Q| \le 2d_Q$, because at most $d_Q$ bases are incident to each endpoint of that edge.
	Otherwise, if there is a pair of non-near bases in $B_Q$, then it follows from Lemma~\ref{lem:near-bases} that they must belong to two gadgets associated with the same $H$-edge.
	Furthermore, by this lemma every other base in $B_Q$ must be near each of these two bases.
	Therefore, by Proposition~\ref{prop:tree-packing} all the other bases in $B_Q$ must form a star which implies that there are at most $d_Q$ such bases.
	Therefore, $|S_Q|=|B_Q| \le 2+d_Q$ in this case.
	If all the gadgets in $S_Q$ are $Q$-pure, then recall that $d_Q \le 1$ and therefore $|S_Q|=|B_Q| \le \max\{2d_Q,2+d_Q\} = 3$. Otherwise, if $d_Q=2$ then $|S_Q| \le 4$.
	
	Let $S_P$ denote the $P$-hooking gadgets in $S$ and let $d_P$ denote the maximum vertex-degree is the subgraph of $U$ consisting of their corresponding base edges in $U$.
	Clearly, the same arguments apply to $S_P$, therefore if all the gadgets in $S$ are pure of the same type, then $|S|\le 3$.
	
	It remains to consider the case that not all gadgets are pure of the same type.
	If there is at least one double gadget, then $|S| \le |S_Q|+|S_P|-1 \le 7$.
	Otherwise, if all the gadgets are pure but not of the same type, then $|S|=|S_Q|+|S_P| \le 3+3 < 7$. This concludes the case that $U$ is a tree.
	
	Suppose now that $\mu(U) > 0$ and let $T$ be a spanning tree of $U$.
	Clearly, $U$ contains $\mu(U)$ edges not in $T$.
	For each such edge choose one of its endpoints and discard from $S$ any gadget that contains this endpoint.
	Denote by $S'$ the set of remaining gadgets.
	Note that, as above, every vertex belongs to at most two gadgets in $S$, thus, at most $2\mu(U)$ gadgets are discarded.
	Furthermore, if all the gadgets are pure, then every vertex belongs to at most one gadget in $S$ and therefore at most $\mu(U)$ gadgets are discarded.
	We would like to apply the above arguments to conclude that $|S'|\le 7$ (and $|S'|\le 3$, if all gadgets are pure of the same type).
	Clearly, $S'$ is admissible.
	Note that $S'$ also satisfies the nearness relations. 
	Indeed, if the deleted edge $(u,v)$ was a witness of two near bases, then one of them is incident to $u$, thus, the corresponding gadget and base were discarded.
	
	It follows from the discussion above that $|S'|\le 7$ and therefore $|S|\le 2\mu(U)+7$.
	If all the gadgets in $S$ are pure of the same type, then every vertex belongs to at most one gadget and therefore at most $\mu(U)$ gadgets are discarded and at most three survive.
	Thus, $|S|\le \mu(U)+3$ in this case.
\end{proof}

\begin{proposition}\label{prop:cycle}
	A connected component $U$ contains a cycle if:\\
	(a)~there are three pure same-type gadgets associated with three $H$-edges connecting $U$ to three distinct components; or\\
	(b)~$U$ has two distinct neighbors $V_1$ and $V_2$ in $H$, there are two gadgets associated with $UV_1$ whose bases in $U$ are vertex-disjoint, there are two gadgets associated with $UV_2$ whose bases in $U$ are vertex-disjoint, and these four gadgets are pure of the same type.
\end{proposition}

\begin{proof}
	Suppose that (a) holds.
	Then by Proposition~\ref{prop:sharing} those gadgets have three vertex-disjoint edge bases in $U$.
	Therefore by Lemma~\ref{lem:near-bases} these base edges are pairwise near and it follows from Proposition~\ref{prop:tree-packing} that $U$ is not a tree.
	
	Suppose now that (b) holds. Considering the four corresponding bases in $U$, it follows from the assumption and Proposition~\ref{prop:sharing} that they are pairwise vertex-disjoint.
	By Lemma~\ref{lem:near-bases} there is an edge of $\GD$ between each of the two bases that corresponds to $UV_1$ and each of the two bases that corresponds to $UV_2$.
	Therefore, by contracting the four bases we obtain a copy of $K_{2,2}$.
	Thus, $\mu(U) \ge 1$.
\end{proof}

\paragraph{Summarizing the properties of $H$.}\label{prop:H}
Before proceeding to the proof of the optimal bound, we summarize the main properties of the transitions graph $H$ for convenience and future reference.
%Furthermore, we view $H$ as a combinatorial object disregarding its relation to the two intersecting polygons.
%This will allow us to apply induction, contract its edges, etc.
Recall that every node of $H$ corresponds to a bipartite graph.
Every edge $UV$ of $H$ is associated with at least one `gadget': two pairs of adjacent vertices (`bases' of the gadget), one from $U$ and the other from $V$. Each such gadget has at least one of the labels $\{Q,P\}$. If it has both labels, then we call it `double' and otherwise it is 'pure'.
The following summarizes Lemma~\ref{lem:same-type-not-disjoint}, Lemma~\ref{lem:near-bases}, Proposition~\ref{prop:sharing} and Proposition~\ref{prop:cut}.

\medskip
$H$ satisfies the following properties:
\begin{enumerate}
	\item There are no two node-disjoint $H$-edges associated with two gadgets sharing a label.
	
	\item Consequently, for each label the $H$-edges associated with gadgets of this label form a star or a triangle.
	
	\item If two gadgets with the same label are associated with two $H$-edges $UV$ and $UZ$ ($V\ne Z$), then their bases in $U$ are near. If both gadgets are pure, then their bases are vertex-disjoint and therefore are connected by an edge of $U$.
	
	\item For every partition of the nodes of $H$ into two disjoint subsets $A$ and $B$ such that $|A| \ge 2$ and $|B| \ge 1$ there are two gadgets crossing the cut with vertex-disjoint bases in $A$.
\end{enumerate}

\section{The optimal bound $X \le mn-(m+n)+3$}
\label{sec:optimal}

Recall that it is enough to show that $c(\GD)-\mu(\GD) \le 4$. 
Therefore, we may assume that $c(\GD)-\mu(\GD) \ge 5$ and in particular $c(\GD) \ge 5$.
We consider the transitions graph $H$ and denote by $H_Q$ (resp., $H_P$) the subgraph of $H$ that consists of $H$-edges associated with $Q$-hooking (resp., $P$-hooking) gadgets.
Recall that it follows from Lemma~\ref{lem:same-type-not-disjoint} that each of these subgraphs is either a triangle or a star (ignoring isolated nodes).
We proceed by considering the possible cases.

\begin{case}
{\bf At least one of $H_Q$ and $H_P$ is a star.}
Suppose without loss of generality that $H_Q$ is a star whose center is $U_Q$.
Then $H_P$ is either a triangle or a star.

\begin{subcase}
{\bf $H_P$ is a triangle.}
%Let $L_1,\ldots,L_k$ be the leaves of the star $H_Q$.
Since $H$ is connected, $H_Q$ and $H_P$ share at least one node.
There are four possible cases: either $U_Q$ is one node of the triangle $H_P$ (see Figure~\ref{fig:star-triangle}(a)), or $U_Q$ is not a node of the triangle $H_P$ and there are exactly $i$ leaves in $H_Q$ that belong to $H_P$, for $i=1,2,3$ (see Figure~\ref{fig:star-triangle}(b)-(d)).
\begin{figure}
	\centering
	\includegraphics[width=16cm]{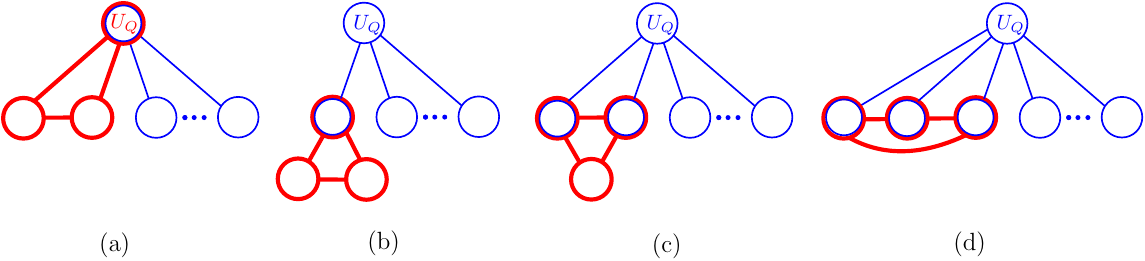}
	\caption{If $H_Q$ is a star and $H_P$ is a triangle then either (a) the center of the star $U_Q$ belongs to $H_P$; or (b)-(d) exactly $i$ of its leaves in $H_Q$ belong to $H_P$, for $i=1,2,3$.}
	\label{fig:star-triangle}
\end{figure}
In the first case, $U_Q$ is connected to $c(\GD)-3$ leaves that do not belong to $H_P$.
For each such leaf $L$ we apply Proposition~\ref{prop:cut} on the cut $(V(H)\setminus \{L\}, \{L\})$ and obtain two $Q$-pure gadgets with disjoint bases in $U_Q$.
This results in an admissible set $S$ of $2(c(\GD)-3)$ $Q$-pure gadgets that satisfies the properties in Proposition~\ref{prop:2-capacities}.
Therefore, $2(c(\GD)-3) = |S| \le \mu(U_Q)+3 \le \mu(\GD)+3 \le c(\GD)-5+3$.
However, this implies that $c(\GD) \le 4$.

Otherwise, suppose that $U_Q$ does not belong to $H_P$ (see Figures~\ref{fig:star-triangle}(b)-(d)).
Then by applying Proposition~\ref{prop:cut} on the cut $(V(H)\setminus \{L\}, \{L\})$ for each leaf $L$ from the $c(\GD)-4$ leaves in $H_Q$ that do not belong to $H_P$ and also for the cut $(V(H)\setminus V(H_P), V(H_P))$, we again obtain an admissible set of $2(c(\GD)-3)$ $Q$-pure gadgets and conclude that $c(\GD) \le 4$ as above.

\begin{comment}
Let $V_1,V_2,V_3$ be the non-isolated vertices of $H_P$.
We could proceed by considering whether $U_Q$ is one of these vertices, but this can be avoided by contracting $V_1,V_2,V_3$ into a single vertex $V$ (which contains the union of the three components).
The resulting graph $H'$ still satisfies the properties has $c(\GD)-2 \ge 3$ vertices, all of its edges are associated with $Q$-pure gadgets and therefore it is a star whose center $U$ is either $U_Q$ or $V$ (the latter holds if $U_Q \in \{V_1,V_2,V_3\}$).
For a leaf $L$ in $H'$ we can apply Proposition~\ref{prop:cut} on the cut $(V(H')\setminus \{L\}, \{L\})$, and obtain two $Q$-pure gadgets with disjoint bases in $U$.
By repeating this for each of the $c(\GD)-3$ leaves in $H'$, we obtain an admissible set $S$ of $2(c(\GD)-3)$ $Q$-pure gadgets that satisfies the properties in Proposition~\ref{prop:2-capacities}.
Therefore, $2(c(\GD)-3) = |S| \le \mu(U)+3 \le \mu(\GD)+3 \le c(\GD)-5+3$.
However, this implies that $c(\GD) \le 4$.
\end{comment}
\end{subcase}

\begin{subcase}
{\bf $H_P$ is a star whose center is $U_Q$.}
In this case $H$ is a star with center $U_Q$ and $c(\GD)-1$ leaves.
For a leaf $L$ we apply Proposition~\ref{prop:cut} on the cut $(V(H)\setminus \{L\}, \{L\})$, and obtain two gadgets with disjoint bases in $U_Q$.
By repeating this for each of the $c(\GD)-1$ leaves, we obtain an admissible set $S$ of $2(c(\GD)-1)$ gadgets that satisfies the properties in Proposition~\ref{prop:2-capacities}.
Note that these gadgets are not necessarily pure.
Therefore, $2(c(\GD)-1) =|S| \le 2\mu(U_Q)+7 \le 2\mu(\GD)+7 \le 2(c(\GD)-5)+7 = 2c(\GD)-3$, which is impossible.
\end{subcase}

\begin{subcase}\label{subcase:diff-centers}
{\bf $H_P$ is a star with center $U_P \ne U_Q$.}
Each of the $c(\GD)-2$ other components is adjacent to $U_Q$, $U_P$, or both.
For each such component $V \ne U_Q,U_P$ we apply Proposition~\ref{prop:cut} on the cut $(V(H)\setminus \{V\},\{V\})$ and obtain two gadgets with disjoint bases in $U_Q \cup U_P$.
Note that each of these gadgets is either $Q$-pure or $P$-pure.
Denote by $n_Q$ and $n_P$ the number of those gadgets, respectively.
Then $n_Q+n_P=2(c(\GD)-2)$.
Applying Proposition~\ref{prop:2-capacities} on $U_Q$ and the subset of those $Q$-pure gadgets, we get that $n_Q \le \mu(U_Q)+3$. 
Similarly, we have $n_P \le \mu(U_P)+3$.
Thus we have, 
\begin{equation}\label{eq:5-components}
2(c(\GD)-2) = n_Q+n_P \le \mu(U_Q)+3 + \mu(U_P)+3 \le \mu(\GD)+6 \le c(\GD)-5+6,
\end{equation}
which implies that $c(\GD) \le 5$.
Since we assumed that $c(\GD) \ge 5$, it remains to consider the case that $\GD$ contains exactly five connected components.

Suppose that $c(\GD)=5$ and let $V_1,V_2,V_3$ be the other connected components of $\GD$ apart from $U_Q$ and $U_P$.
We may further assume that $\mu(U_Q)=\mu(U_P)=\mu(\GD)=0$ and that~(\ref{eq:5-components}) holds with equality, for otherwise we immediately get $c(\GD)-\mu(\GD)<5$ and we are done. 
Thus, $n_Q+n_P=6$ and therefore $n_Q=n_P=3$ (because $n_x \le \mu(U_x) + 3 = 3$, for $x \in \{P,Q\}$).

Recall that for every $i=1,2,3$, exactly two of those six gadgets have one base in $V_i$, the other base in $U_Q$ or $U_P$, and if both gadgets are associated with the same $H$-edge, then their bases in the star-center are disjoint.
Observe first that it follows from the first part of Proposition~\ref{prop:cycle} that it is impossible to have three gadgets associated with the three edges $U_QV_1,U_QV_2,U_QV_3$, since this would imply that $\mu(U_Q) \ge 1$.
By the second part of Proposition~\ref{prop:cycle} it is also impossible that two of these gadgets are over $U_QV_i$ and two others are over $U_QV_j$, for $i\ne j$.
Applying the same arguments for $U_P$, we conclude that, without loss of generality, there are two gadgets over each of $U_QV_1$ and $U_PV_3$ and one gadget over each of $U_QV_2$ and $U_PV_2$ (see Figure~\ref{fig:5-path}).

By applying Proposition~\ref{prop:cut} on the cut $(\{U_Q,V_1,V_2\},\{U_P,V_3\})$, it follows that there must be another $H$-edge crossing this cut except $U_PV_2$.
Such an edge can only be $U_QU_P$ (see Figure~\ref{fig:5-path}).
\begin{figure}
	\centering
	\includegraphics[width=12cm]{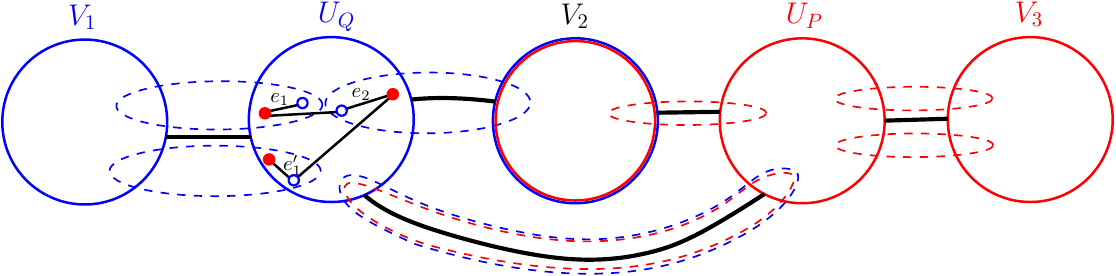}
	\caption{An illustration for Subcase~\ref{subcase:diff-centers}: $c(\GD)=5$, $U_Q$ is connected by two $Q$-pure gadgets to $V_1$, by a $Q$-pure gadget to $V_2$, and by a double gadget to $U_P$. $U_P$ is also connected by two $P$-pure gadgets to $V_3$ and by a $P$-pure gadget to $V_2$.}
	\label{fig:5-path}
\end{figure}
Note that no gadget associated with this edge can be pure. 
Indeed, if it is $Q$-pure, then considering pure gadgets between $U_Q$ and each of $V_1$, $V_2$ and $U_P$, it follows that $\mu(U_Q)\ge 1$ by the first part of Proposition~\ref{prop:cycle}.
Similarly, $U_QU_P$ cannot be associated with a $P$-pure gadget.
Let $e_1$ and $e'_1$ be vertex-disjoint bases in $U_Q$ from $Q$-pure gadgets associated with $U_QV_1$, and let $e_2$ be a base in $U_Q$ from a $Q$-pure gadget associated with $U_QV_2$.
Note that $e_2$ does not share vertices with $e_1$ and $e'_1$ by Proposition~\ref{prop:sharing}, and that there are edges connecting $e_2$ to both of $e_1$ and $e'_1$ by Lemma~\ref{lem:near-bases}.
Therefore $e_1$ and $e'_1$ are not connected by an edge, for otherwise we could contract $e_1$, $e'_1$ and $e_2$ and get a $3$-cycle implying $\mu(U_Q) \ge 1$.

Consider now the cut $(\{U_Q,V_1\},\{U_P,V_2,V_3\})$.
Every gadget that crosses this cut is either a $Q$-pure gadget associated with $U_QV_2$ or a double gadget associated with $U_QU_P$.
By Lemma~\ref{lem:near-bases} the base edge in $U_Q$ of such gadgets must be near each of $e_1$ and $e'_1$.
It follows from Proposition~\ref{prop:tree-packing}(b) that these base edges form a star. 
However, this contradicts Proposition~\ref{prop:cut} by which there should be two vertex-disjoint base edges in $U_Q$ that belong to two gadgets crossing this cut.
\end{subcase}
\end{case}

\begin{case}\label{subcase:2triangles}
{\bf Both $H_Q$ and $H_P$ are triangles.}
Since $H=H_Q \cup H_P$ is connected, it follows that $c(\GD) \le 5$.
Recalling the assumption that $c(\GD) \ge 5$, we may assume that the two triangles share exactly one node, denote it by $U$.
Let $V_1$ and $V_2$ be the other nodes of $H_Q$ and let $V_3$ and $V_4$ be the other nodes of $H_P$.
Note that since no $H$-edge belongs to both $H_Q$ and $H_P$, it follows that every gadget associated with the six edges of $H$ is pure.
We may also assume that each of the five connected components of $\GD$ is a tree, for otherwise $\mu(\GD) \ge 1$ and we are done.

Let $S_i$ be the set of endpoints of base edges in $U$ of gadgets that are associated with $UV_i$, for $i=1,2,3,4$.
It follows from Proposition~\ref{prop:sharing} that $S_i \cap S_j = \emptyset$ for $i\ne j$.
By Lemma~\ref{lem:near-bases}, each of the subgraphs of $U$ induced by $S_1\cup S_2$ and $S_3 \cup S_4$ is connected.
Recall that $U$ is a tree, hence there is a unique path connecting the subtrees $U[S_1\cup S_2]$ and $U[S_3 \cup S_4]$.
Let $q^*$ be a $Q$-side that corresponds to a vertex on this path. See Figure~\ref{fig:2triangles} for an illustration.
\begin{figure}
	\centering
	\includegraphics[width=8cm]{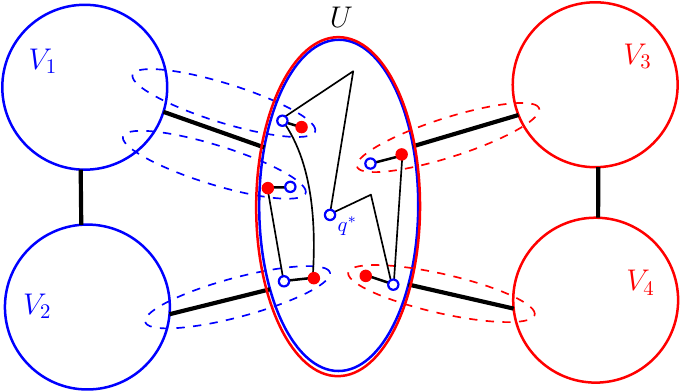}
	\caption{An illustration for Case~\ref{subcase:2triangles}: $H_Q$ and $H_P$ are two triangles sharing exactly one node $U$.}
	\label{fig:2triangles}
\end{figure}

Omitting $q^*$ from $Q$, we obtain a sequence of $n-1$ $Q$-sides.
Some of them belong to $V_1 \cup V_2$ and some to $V_3 \cup V_4$, as every connected component of $\GD$ contains sides of both polygons and we omitted a side from $U$.
Recall that transitions along this sequence between sides that belong to different connected components correspond to $H$-edges.
Since there are no transitions between $\{V_1,V_2\}$ and $\{V_3,V_4\}$, it follows that there is a subsequence of consecutive $Q$-sides $q'_1,\ldots,q'_k$ such that each of them belongs to $U$, $q'_1 \in S_1\cup S_2$, $q'_k \in S_3 \cup S_4$, the side before $q'_1$ belongs to $V_1 \cup V_2$ and the side after $q'_k$ belongs to $V_3 \cup V_4$ (we may assume that this is the order without loss of generality, since we can reverse the sequence of sides if needed). 

We claim that for every $i=1,\ldots,k-1$, there is a side of $P$ that is disjoint from both $q'_i$ and $q'_{i+1}$.
Indeed, suppose that there is no $P$-side that is disjoint from both $q'_i$ and $q'_{i+1}$.
Then, it follows from Lemma~\ref{lem:gadget} that $q'_i$ and $q'_{i+1}$ are part of a gadget $g$. 
Because $q'_i,q'_{i+1} \in U$, all the sides of $g$ belong to $U$.
Since there is a $Q$-pure gadget associated with $V_1V_2$ and its sides are disjoint from those of $g$, it follows from Lemma~\ref{lem:eight} that $g$ is not $Q$-hooking.
Similarly, since there is a $P$-pure gadget associated with $V_3V_4$ it follows that $g$ is also not $P$-hooking, a contradiction.

Therefore, for every $i=1,\ldots,k-1$, there is a $P$-side $p'_i$ (in $U$) such that $(q'_i,p'_i), (q'_{i+1},p'_i) \in E(\GD)$.
Thus, there is a path between $S_1 \cup S_2$ and $S_3 \cup S_4$, which avoids $q^*$, a contradiction.
Hence, $U$ is not a tree and we are done.
\end{case}

\bibliographystyle{plainurl}% the mandatory bibstyle
\bibliography{b}

\appendix
\section{The gap in the proof of the upper bound $X \le mn-(m+n)+C$ in~\cite{AKR22}}
\label{app:gap}

\begin{figure}[t]
	\centering
	\subfloat[$b_3$ is implicitly assumed to cross $A_3$ between $A_3 \cap a_1$ and $A_3 \cap a_4$.]{\includegraphics[width= 7cm]{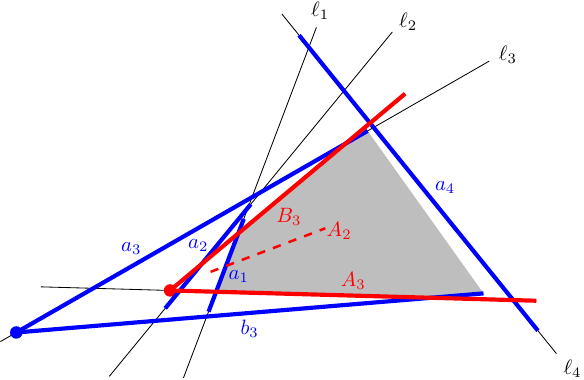}\label{fig:gap1}}
	\hspace{5mm}
	\subfloat[the case that $b_3$ crosses $A_3$ between $A_3 \cap a_2$ and $A_3 \cap a_1$ is not handled.]{\includegraphics[width= 7cm]{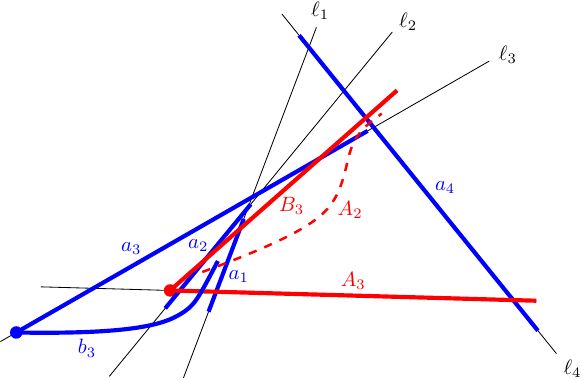}\label{fig:gap2}}
	\caption{The gap in the proof of the upper bound in~\cite{AKR22}.}
	\label{fig:gap}
\end{figure}

\begin{figure}[t]
	\centering
	\includegraphics[width= 16cm]{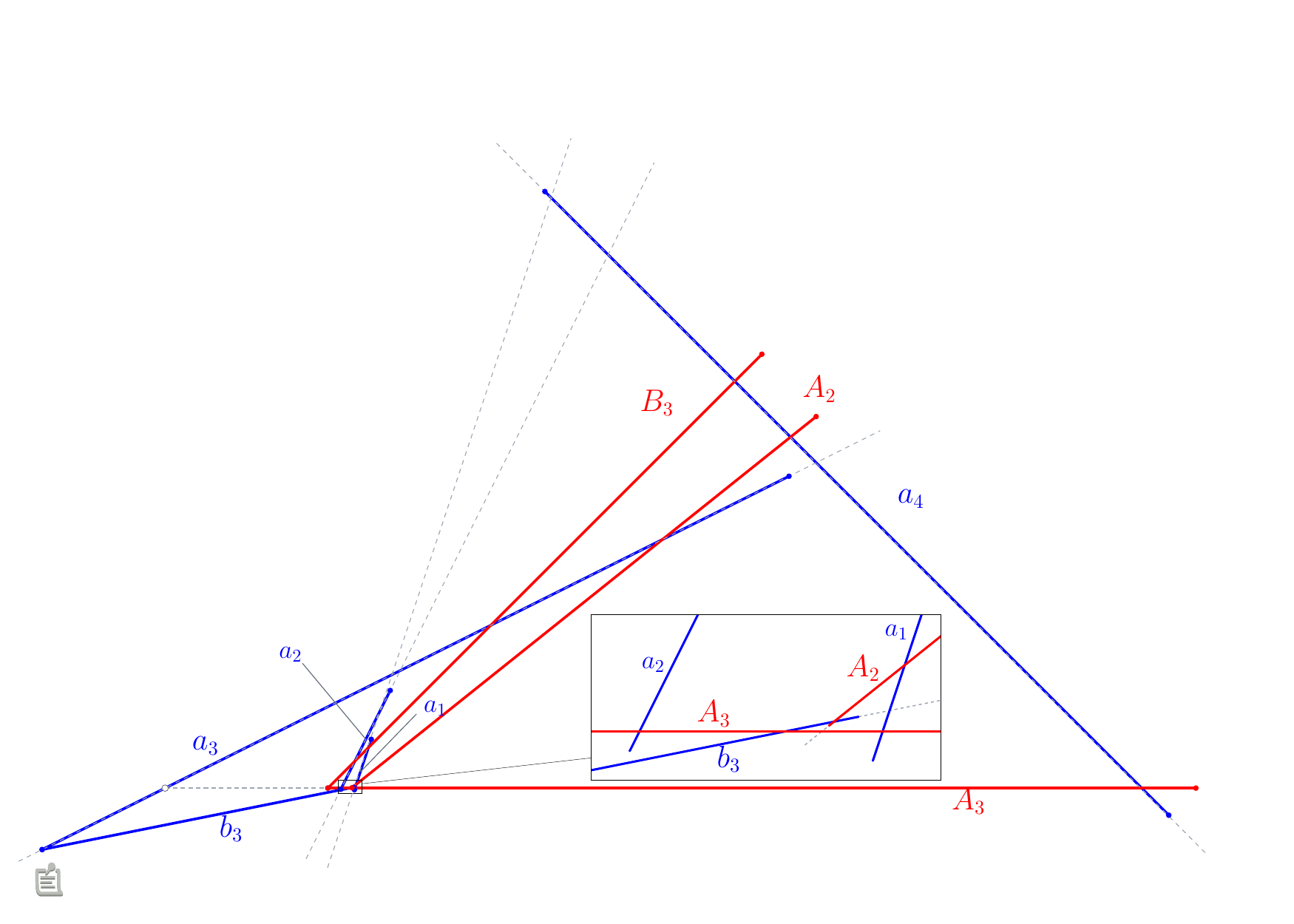}
	\caption{A realization of the omitted case in~\cite{AKR22} by straight-line segments.}
	\label{fig:gap3}
\end{figure}

For the proof of the upper-bound in~\cite{AKR22} one assumes that there are more than $C$ connected components in $\GD$, where $C<2^{2^{70}}$ is a huge constant.
Using Ramsey theory and Lemma~\ref{lem:eight} it then follows that there are six pairwise side-disjoint $Q$-hooking gadgets $(A_i,B_i,a_i,b_i)$, $i=1,2,\ldots,6$, with the following properties (some are implied by the other):
\begin{itemize}
	\item The lines $\ell_i$ containing the sides $a_i$ form a $6$-cap such that the slope of $\ell_i$ is greater than the slope of $\ell_j$ for every $i<j$;
	\item the sides $a_i$ are pairwise disjoint and $\ell_i \cap \ell_j$ lies on $a_j$, for every $i < j$;
	\item $a_i$ and $b_i$ belong to different connected components of $\GD$ for every $i$;
	\item $a_i$ and $a_j$ belong to different connected components of $\GD$ for every $i\ne j$;
	\item $b_i$ and $b_j$ belong to the same connected component of $\GD$ for every $i$ and $j$;
	\item $A_i$ crosses $a_j$ for every $i\ne j$;
	\item $A_i$ crosses $b_j$ for every $i$ and $j$; and
	\item $B_i$ crosses $a_j$ for every $i$ and $j$.
\end{itemize}
The proof proceeds by a case-analysis resulting in a contradiction in every case.
In Case~2, one assumes that for $1 \le i \le 4$ the side $A_i$ does not cross $\ell_i$ to the left of $a_i$.
Then one considers the four possible subcases depending on whether $a_3 \cap b_3$ is the left or right endpoint of $a_3$ and whether $A_3 \cap B_3$ is the left or right endpoint of $A_3$.
In Subcase~2.2 $a_3 \cap b_3$ is the left endpoint of $a_3$ and $A_3 \cap B_3$ is the left endpoint of $A_3$.
The side $a_1$ crosses both $A_3$ and $B_3$, which partition it into three parts.
It is then claimed that if $A_2$ intersects the middle part of $a_1$ (between $a_1 \cap A_3$ and $a_1 \cap B_3$), then $A_2$ cannot cross both $a_3$ and $b_3$ in the wedge determined by $a_3$ and $b_3$.
That is true if $b_3$ crosses $A_3$ between $A_3 \cap a_1$ and $A_3 \cap a_4$, as is assumed implicitly by the accompanying figure (see Figure~\ref{fig:gap1}).
However, it is possible that $b_3$ crosses $A_3$ between $A_3 \cap a_2$ and $A_3 \cap a_1$, and then $A_2$ may cross both $a_3$ and $b_3$ inside the above-mentioned wedge, see Figures~\ref{fig:gap2} and~\ref{fig:gap3}.

It seems that the gap in this case-analysis can be resolved by considering seven gadgets instead of six (thus taking a greater constant $C$). However, since the current proof is simpler and gives a better upper bound, there is not much point in going into those technicalities. 
\end{document}